\documentclass[a4paper,9pt]{article}
\usepackage[cp1251]{inputenc}
\usepackage[english]{babel}
\usepackage{amsthm,amsmath,amssymb,latexsym}
\theoremstyle{definition}

\theoremstyle{remark}

\theoremstyle{theorem}

\theoremstyle{definition}

\theoremstyle{remark}

\newcommand{\ZZ}{{\mathbb Z}}

\newcommand{\Ann}{\mathop{\rm Ann}\nolimits}

\newcommand{\ov}{\overline}

\begin{document}

\begin{center}
\large{\textbf{Cayley-Dickson Process and\\ Centrally Essential Rings}}
\end{center}

\hfill {\sf V.T. Markov}

\hfill Lomonosov Moscow State University

\hfill e-mail: vtmarkov@yandex.ru

\hfill {\sf A.A. Tuganbaev}

\hfill National Research University "MPEI"

\hfill Lomonosov Moscow State University

\hfill e-mail: tuganbaev@gmail.com

{\bf Abstract.}
We describe associative center $N(R)$ and the center $Z(R)$ of the ring $R$ obtained
by applying the generalized Cayley-Dickson construction and we find conditions
under which the ring $R$ is $N$-essential or centrally essential.
The obtained results are applied to generalized quaternion rings and octonion rings; we use them to construct an example of a non-associative centrally essential ring.

V.T.Markov is supported by the Russian Foundation for Basic Research, project 17-01-00895-A. A.A.Tuganbaev is supported by Russian Scientific Foundation, project 16-11-10013.


\section{Introduction}\label{section1}
All considered rings are unital but not necessarily associative.

We denote by $N(R)$, $K(R)$, $Z(R)$ the associative center, the commutative center and the center (in the sense of \cite[$\S$7.1]{Zhevlakov}) of the ring $R$,  respectively; the definitions are also given at the end of the introduction. It is clear that $N(R)$, $Z(R)$ are subrings in $R$ and the ring $R$ is a unitary (left and right) $N(R)$-module and $Z(R)$-module.

We denote by $[A,A]$ the ideal of the ring $A$ generated by commutators of all its elements.

\textbf{1.1. The definitions of a centrally essential ring and an $N$-essential ring.}\\ A ring $R$ is said to be \emph{centrally essential}\label{Z-ess} if $Z(R)r\cap Z(R)\neq 0$ for any non-zero element $r\in R$, i.e., $Z=Z(R)$ is an essential submodule of the module ${}_ZR$.\\
A ring $R$ is said to be \emph{left $N$-essential}\label{N-ess} if
$N(R)r\cap N(R)\neq 0$ for any non-zero element $r\in R$, i.e., $N=N(R)$ is an
essential submodule of the module ${}_NR$.

The following definition several generalizes the definition of the Cayley-Dickson process given in \cite[\S 2.2]{Zhevlakov}, see \cite{Albert}.

\textbf{1.2. The Cayley-Dickson process.}\label{pCD} Let $A$ be a ring with involution $*$ and $\alpha$ an invertible symmetrical element of the center of the ring $A$.
We define a multiplication operation on the Abelian group $A\oplus A$ as follows:
\begin{equation}\label{mCD}
(a_1,a_2)(a_3,a_4) = (a_1a_3 +\alpha
a_4a_2^*\;,\; a_1^*a_4 + a_3a_2)
\end{equation}
for any $a_1,\ldots,a_4\in A$. We denote the obtained ring by $(A,\alpha)$.

The elements of the ring $(A,\alpha)$ of the form $(a,0)$, $a\in A$, form a
subring in ring $(A,\alpha)$ which is isomorphic to the ring $A$; we will identify
them with the corresponding elements of the ring $A$. We set
$\nu=(0,1)\in (A,\alpha)$. Then $a*\nu=(0,a)=\nu a$ for any
$a\in A$ and $\nu^2=\alpha$. Thus, $(A,\alpha)=A+A\nu$.

The following properties are directly verified with the use of relation \eqref{mCD}.\label{elprop}

\textbf{1.2.1.} $\nu^2=\alpha$ and $\nu a=a^*\nu$ for any element $a\in A$.

\textbf{1.2.2.} $(1,0)$ is the identity element of the ring $(A,\alpha)$.

\textbf{1.2.3.} The set $\{(a,0)\,|\, a\in A\}$ is a subring of the ring $(A,\alpha)$ which is isomorphic to the ring $A$.

\textbf{1.2.4.} The mapping $(a,b)\mapsto (a^*,-b)$, $a,b\in A$, is an involution of the ring $(A,\alpha)$.

The main results of this paper are Theorems 1.3, 1.4 and 1.5.

\textbf{Theorem 1.3.}\label{main}
Let $A$ be a ring with center $C=Z(A)$, $I=\Ann_C([A,A])$, $R=(A,\alpha)$.\\
\textbf{1.3.1.} $N(R)=\{(x,y)\colon x\in C,\;y\in I\}$.\\
\textbf{1.3.2.} A ring $R$ is left (right) $N$-essential if and only if
$A$ is a centrally essential ring and $I$ is an essential ideal of the ring $C$.

\textbf{Theorem 1.4.}\label{main2}
Let $A$ be a ring with center $C=Z(A)$, $I=\Ann_C([A,A])$, $B=\{a\in C\colon a=a^*\}$, $J=\Ann_B(\{a-a^*\;:\;a\in A\})$.\\
\textbf{1.4.1.} $Z(R)=\{(x,y)\;|\;x\in B\cap C,\;y\in I\cap J\}$.\\
\textbf{1.4.2.} A ring $R$ is a centrally essential if and only if $B$ is an
essential $B$-submodule of the ring $R$ and $J$ is an essential ideal of the ring $B$.

\textbf{Theorem 1.5.}\label{main3}\\
\textbf{1.5.1.} There exists a finite non-associative and non-commutative alternative centrally essential
ring.\\
\textbf{1.5.2.} There exists a finite non-commutative and non-alternative centrally essential ring.

The proof of Theorems 1.3, 1.4 and 1.5 is decomposed into several assertions given
in consequent sections. Some of these assertions are of independent interest.

We give some necessary definitions and notation. The \emph{associator} of three elements $a,b,c$ of the ring $R$ is the element $(a,b,c)=(ab)c-a(bc)$ and the \emph{commutator} of two elements $a,b\in R$ is the element $[a,b]=ab-ba$. For ring $R$, the \emph{associative center, commutative center and center} of $R$ are the sets
$$
\begin{array}{l}
N(R)=\{x\in R:\forall a,b\in R,\;(x,a,b)=(a,x,b)=(a,b,x)=0\},\\
K(R)=\{x\in R:\forall a\in R,\;[x,a]=0\},\\
Z(R)=N(R)\cap K(R),
\end{array}
$$
respectively

If $M$ is a left module over the ring $R$ and $S$ is a subset of the module $M$, then
$\Ann_R(S)$ denotes the annihilator of the set $S$ in the ring $R$, i.e.,
$\Ann_R(S)=\{r\in R: rS=0\}$.

A ring $R$ is said to be \emph{alternative} if $(a,a,b)=(a,b,b)$ for any
$a,b\in R$. By the Artin theorem \cite[Theorem 2.3.2]{Zhevlakov}, the ring $R$
is alternative if and only if any two elements of $R$ generate the associative subring.

\section{Associative Center of the Ring Obtained by Applying Cayley-Dickson Process}
\label{section2}
In Section 2, we fix a ring $A$ and an element $\alpha$ which satisfy the conditions of the Cayley-Dickson process from 1.2 \label{pCD}; we also set $R=(A,\alpha)$.

\textbf{Lemma 2.1.}\label{ident}
An element $(x,y)\in R$ belongs to the ring $N(R)$ if and only if for any two elements $u,v\in A$, the following two relation systems  hold
\begin{equation}\label{id_x}
\begin{array}{l}
(xu)v=x(uv),(ux)v=u(xv),(uv)x=u(vx),\\
v(ux)=x(vu),(xu)v=u(vx),(vu)x=(xv)u,\\
v(xu)=(vu)x,v(ux)=(vx)u,x(uv)=u(xv),\\
(ux)v=(uv)x,v(xu)=(xv)u,x(vu)=(vx)u;
\end{array}
\end{equation}
\begin{equation}\label{id_y}
\begin{array}{l}
(uy)v=y(vu),(uy)v=(yv)u,y(vu)=u(yv),\\
v(yu)=y(uv),(yu)v=(vy)u,y(uv)=(vy)u,\\
v(uy)=(uv)y,v(uy)=u(vy),(vu)y=u(vy),\\
(yu)v=(vu)y,v(yu)=u(yv),(uv)y=(yv)u.
\end{array}
\end{equation}

%

\textbf{Proof.} Let $(x,y)\in R$. Since associators are linear, \mbox{$(x,y)\in N(R)$} if and only if for any two elements $u,v\in A$, we have
\begin{equation}\label{associators}
\begin{array}{l}
((x,y)(u,0)(v,0))=((u,0),(x,y),(v,0))=((u,0),(v,0),(x,y))=0,\\
((x,y)(u,0)(0,v))=((u,0),(x,y),(0,v))=((u,0),(0,v),(x,y))=0,\\
((x,y)(0,u)(v,0))=((0,u),(x,y),(v,0))=((0,u),(v,0),(x,y))=0,\\
((x,y)(0,u)(0,v))=((0,u),(x,y),(0,v))=((0,u),(0,v),(x,y))=0.
\end{array}
\end{equation}
By calculating associators from \eqref{associators}, we obtain the following system consisting of 12 relations
$$
\begin{array}{l}
((xu)v,v(uy))=(x(uv),(uv)y),\\
((ux)v,v(u^*y))=(u(xv),u^*(vy)),\\
((uv)x,(v^*u^*)y)=(u(vx),u^*(v^*y)),\\
({\alpha}v(y^*u^*),(u^*x^*)v)=({\alpha}(u^*v)y^*,x^*(u^*v)),\\
({\alpha}v(y^*u),(x^*u^*)v)=({\alpha}u(vy^*),u^*(x^*v)),\\
({\alpha}y(v^*u),x(u^*v))=({\alpha}u(yv^*),u^*(xv)),\\
({\alpha}(uy^*)v,v(x^*u))=({\alpha}(vu)y^*,x^*(vu)),\\
({\alpha}(yu^*)v,v(xu))=({\alpha}(vy)u^*,(xv)u),\\
({\alpha}y(u^*v^*),x(vu))=({\alpha}(v^*y)u^*,(vx)u),\\
({\alpha}v(u^*x),{\alpha}(yu^*)v)=({\alpha}x(vu^*),{\alpha}(vu^*)y),\\
({\alpha}v(u^*x^*),{\alpha}(uy^*)v)=({\alpha}(x^*v)u^*,{\alpha}(vy^*)u),\\
({\alpha}(vu^*)x,{\alpha}(uv^*)y)=({\alpha}(xv)u^*,{\alpha}(yv^*)u).
\end{array}
$$
By equating components of equal elements of the ring $R$ and considering that the element $\alpha$ is invertible, we obtain the following equivalent system
$$
\begin{array}{l}
(xu)v=x(uv),v(uy)=(uv)y),
(ux)v=u(xv),v(u^*y)=u^*(vy),\\
(uv)x=u(vx),(v^*u^*)y=u^*(v^*y),
 v(y^*u^*)=(u^*v)y^*,(u^*x^*)v=x^*(u^*v),\\
 v(y^*u)=u(vy^*),(x^*u^*)v=u^*(x^*v),
 y(v^*u)=u(yv^*),u^*(xv)=x(u^*v),\\
 (uy^*)v=(vu)y^*,v(x^*u)=x^*(vu),
 (yu^*)v=(vy)u^*,v(xu)=(xv)u,\\
 y(u^*v^*)=(v^*y)u^*,x(vu)=(vx)u,
 v(u^*x)=x(vu^*),(yu^*)v=(vu^*)y,\\
 v(u^*x^*)=(x^*v)u^*,(uy^*)v=(vy^*)u,
 (vu^*)x=(xv)u^*,(uv^*)y=(yv^*)u.
\end{array}
$$
We replace the equations, both parts of which contain $x^*$ or $y^*$, by relations of
conjugate elements. We note that either $u$ or $u^*$ stands in every equation. Therefore, we can put $u$ instead of $u^*$, since $A^*=A$. Similarly, we replace $v^*$ by $v$. By choosing equations containing $x$, we obtain \eqref{id_x},
the remaining equations form the system \eqref{id_y}.~\hfill$\square$

\textbf{Lemma 2.2.}\label{x_in_Z}
Let $x\in A$. The relations \eqref{id_x} hold for all $u,v\in A$ if and only if $x\in Z(A)$.

\textbf{Proof.} Let $x\in A$ and let relations~\eqref{id_x} hold for all $u,v\in A$. The first three relations mean that $x\in N(A)$. It follows from the fourth relation for $u=1$ we obtain $x\in K(A)$. Consequently, $x\in Z(A)$.

Conversely, if $x\in Z(A)$, then each of the relations in \eqref{id_x} is transformed into one of the true relations $x(uv)=x(uv)$ or $x(vu)=x(vu)$, i.e., relations
\eqref{id_x} hold for all $u,v\in A$.~\hfill$\square$

\textbf{Lemma 2.3.}\label{y_in_Ann_I} Let $y\in A$. The relations \eqref{id_y}
hold for all $u,v\in A$ if and only if $y\in \Ann_{Z(A)}([A,A])$.

\textbf{Proof.} Let $y\in A$ and let relations \eqref{id_y} hold for all $u,v\in A$. First of all, we note that for $v=1$, the first equation of \eqref{id_y}
turns into equation $uy=yu$; this is equivalent to the inclusion $y\in K(A)$, since the element$u$ of $A$ is arbitrary.

We verify that $y\in N(A)$. For any two elements $u,v\in A$, we have
$$
\begin{array}{l}
(yu)v\stackrel{\ov 1}{=}(vy)u\stackrel{\ov 2}{=}y(uv),\\
(uy)v\stackrel{1}{=}y(vu)\stackrel{3}{=}u(yv)\\
(uv)y=y(uv)\stackrel{4}{=}v(yu)=v(uy)\stackrel{8}{=}u(vy).
\end{array}
$$
In these transformations, the number over the relation sign is the number of used 
equation in \eqref{id_y} (equations are numbered in rows from the left to right beginning with the first row). The  number underscore denotes that, instead of the given equation, we use the equivalent equation obtained by the permutation of the variables $u,v$.

Consequently, $y\in N(A)\cap K(A)=Z(A)$.

Finally, we take into account the proven to see that already the first equation of \eqref{id_y} implies that $y[u,v]=0$ for any $u,v\in A$, i.e., $y\in\Ann_C([A,A])$.

Conversely, if $y\in \Ann_{Z(A)}([A,A])$, then each of the relations \eqref{id_y} is transformed into the true relation $y(uv)=y(vu)$, i.e., relations
\eqref{id_y} hold for all $u,v\in A$.~\hfill$\square$

\textbf{Remark 2.4.} Theorem 1.3.1\label{main} follows from Lemma 2.1, Lemma 2.2 and Lemma 2.3.\label{ident}\label{x_in_Z}\label{y_in_Ann_I}.

\textbf{Remark 2.5.} It follows from the above classical result (cf.
\cite[Exercise 2.2.2(a)]{Zhevlakov}):\label{crit}\\
A ring $R=(A,\alpha)$ is associative if and only if the ring $A$ is associative and commutative.

\section{$N$-Essentiality Criterion of the Ring Obtained by the Cayley-Dickson Process}

\textbf{Lemma 3.1.}\label{ess_in_R} Let $B$ be a subring of the center of the ring $A$ and $I$ an essential ideal of the ring $B$. If $B$ is an essential $B$-submodule of the module ${}_BA$, then $I$ is an essential $B$-submodule of the module ${}_BR$.

\textbf{Proof.} If $r$ is a non-zero element of the ring $R$, then there exists an element $b\in B$ with $0\neq br\in B$. Therefore, there exists an element $d\in B$ such that $0\neq dcr\in I$ and $Br\cap I\neq 0$.~\hfill$\square$

\textbf{3.2. The proof of Theorem 1.3.2.} Let ring $A$ and the element $\alpha$ satisfy the conditions of 1.2 (the Cayley-Dickson process).\label{pCD} We set $C=Z(A)$ and $I=\Ann_C([A,A])$. It is obvious that $C^*=C$, $I^*=I$, $\alpha C=C$ and $\alpha I=I$.

Let the ring $R=(A,\alpha)$ be $N$-essential. Then for any non-zero element $a\in A$ there exists an element $(x,y)\in N(R)$ such that $(x,y)(a,0)=(xa,ay)\in N(R)\setminus\{0\}$. By Theorem 1.3.1\label{main}, $x\in C$ and $y\in I$. If $xa\neq 0$, then $xa\in C\setminus\{0\}$; otherwise, $ya\in C\setminus\{0\}$. In the both cases, $Ca\cap C\neq 0$. Thus, $A$ is a centrally essential ring.

We prove that $I$ is an essential ideal of the ring $C$. Let $c\in C\setminus \{0\}$. If $Ic\neq 0$, then $Ic\subseteq I$ and $Cc\cap I\neq 0$. Let $Ic=0$. We consider element $(0,c)$. There exists an element $(x,y)\in N(R)$ such that $(x,y)(0,c)=(\alpha cy,x^*c)\in N(R)\setminus\{0\}$. Since $\alpha y\in I$, $\alpha c y=0$, we have that $x^*c\neq 0$ and $x^*c\in I$ by Remark 2.4. Consequently, $Cc\cap I\neq 0$, which is required.

Conversely, let's assume that $A$ is a centrally essential ring and $I$ is a essential ideal in $C$.

Let $(x,y)\in R\setminus\{0\}$. There exists an element $c\in C$ such that $cx\in C\setminus\{0\}$. Since $(c,0)\in N(R)$, we have $0\neq (c,0)(x,y)=(cx,c^*y)\in N(R)(x,y)$. If $c^*y=0$, then $0\neq (cx,0)\in N(R)(x,y)\cap N(R)$. If $c^*y\neq 0$, then by Lemma 3.1\label{ess_in_R} (for $B=C$) there exists an element $d\in C$ such that $dc^*y\in I\setminus \{0\}$. Then
$$
(d^*,0)(c,0)(x,y)=(d^*,0)(cx,c^*y)=(d^*cx,dc^*y)\in N(R)(x,y)\cap N(R)\setminus\{0\}.
$$
Thus, ring $R$ is a $N$-essential.~\hfill$\square$

\section{Proof of Theorem 1.4}\label{main2}\label{section4}
We fix the ring $A$ with center $C=Z(A)$ and element $\alpha$ which satisfy 1.2 (the Cayley-Dickson process).\label{pCD} Let $R=(A,\alpha)$,
\begin{equation}\label{centerdefs}
\begin{array}{l}I=\Ann_C([A,A]),\quad
B=\{a\in C:\;a=a^*\},\\
J=\Ann_B(\{a-a^*\;:\;a\in A\}).
\end{array}
\end{equation}
We note that the sets $B$ and $J$ are invariant with respect to the involution and are closed with respect to the multiplication by $\alpha$.

\textbf{Proposition 4.1.}\label{main2.1}
$Z(R)=\{(x,y)\,|\, x\in B,\,y\in I\cap J\}$.

\textbf{Proof.} Let $(x,y)\in Z(R)$. Since $Z(R)\subseteq N(R)$, it follows from
Theorem 1.3 \label{main} that $x\in C$ and $y\in I$. The relations $(0,1)(x,y)=(x,y)(0,1)$ imply the relations $\alpha y=\alpha y^*$ and $x=x^*$. Consequently, $x\in B$ and $y\in B\cap I$. Next, the relation $(a,0)(x,y)=(x,y)(a,0)$,  $a\in A$, implies the relations $ax=xa$ and $ay=a^*y$. The first relation holds for any $x\in C$ and the second relation means that $y(a-a^*)=0$, i.e., $y\in J$. Consequently, $y\in I\cap J$.

Conversely, if $x\in B$ and $y\in I\cap J$, then $(x,y)\in N(R)$ and for any $a,b\in A$
we have
\begin{equation*}
\begin{array}{l}
(x,y)(a,b)=(xa+\alpha by^*,
x^*b+ay)=(xa+\alpha yb, xb+ay),\\
(a,b)(x,y)=(ax+\alpha
yb^*,a^*y+xb)=(ax+\alpha yb,xb+ay)
\end{array}
\end{equation*}
Thus, $(x,y)\in K(R)$, whence $(x,y)\in Z(R)$.~\hfill$\square$

\textbf{Proposition 4.2.}\label{main2.2}
A ring $R=(A,\alpha)$ is centrally essential if and only if $B$ is a essential $B$-submodule of the ring $R$ and $J'=J\cap I$ is an essential ideal of the ring $B$.~\hfill$\square$

\textbf{Proof.} Let the ring $R=(A,\alpha)$ be centrally essential. Then for any $a\in A\setminus\{0\}$, there exists an element $(x,y)\in Z(R)$ such that $(x,y)(a,0)=(xa,ay)\in Z(R)\setminus\{0\}$. By Proposition 4.1\label{2.1}, $x,xa\in B$ and $y,ay=ya\in J'$. If $xa\neq 0$, then $xa\in B\setminus\{0\}$; otherwise, $ya\in B\setminus\{0\}$. In the both cases, we have $Ba\cap B\neq 0$. Thus, $B$ is a essential submodule of the module ${}_BA$.

We prove that $J'=$ is an essential ideal of the ring $B$. Let $b\in B\setminus \{0\}$. If $J'b\neq 0$, then $J'b\subseteq J'$ and $Bb\cap J'\supseteq J'b\cap J'\neq
0$. Let $J'b=0$. We consider element $(0,b)$. There exists an element $(x,y)\in Z(R)$ such that $(x,y)(0,b)=(\alpha by,x^*b)\in Z(R)\setminus\{0\}$. Since $\alpha
y\in J'$, $\alpha b y=0$, we have that $x^*b\neq 0$, $x\in B$ and $x^*b=xb\in J'$, by Proposition 4.1\label{main2.1}. Consequently, $Bb\cap J'\neq 0$, which is required.

Conversely, let's assume that $B$ is an essential $B$-submodule of the ring $R$ and $J'$ is a essential ideal of the ring $B$.

Let $(x,y)\in R\setminus\{0\}$. First, we assume that $x\neq 0$. There exists an element $b\in B$ such that $bx\in B\setminus\{0\}$. Since $(b,0)\in Z(R)$, we have $0\neq
(b,0)(x,y)=(bx,b^*y)\in Z(R)(x,y)$. If $b^*y=0$, then $0\neq (bx,0)\in Z(R)(x,y)\cap Z(R)$. If $b^*y\neq 0$, then by Lemma 3.1\label{ess_in_R} (for $I=J'$) there exists an element $d\in B$ such that $db^*y\in J'\setminus \{0\}$. Then
$$
(d^*,0)(b,0)(x,y)=(d^*,0)(bx,b^*y)=(d^*bx,db^*y)\in Z(R)(x,y)\cap Z(R)\setminus\{0\}.
$$
Now let $x=0$. Then $y\neq 0$, and there exists an element $d\in B$ such that $dy\in J'\setminus \{0\}$. We obtain $(d^*,0)(0,b)=(0,db)\in Z(R)\setminus \{0\}$ and $(d^*,0)\in Z(R).$ Thus, the ring $R$ is centrally essential.~\hfill$\square$

\textbf{4.3. The completion of the proof of Theorem 1.4.\label{main2}}\\
The first assertion of Theorem 1.4\label{main2} follows from Proposition 4.1\label{main2.1}.\\The second assertion of Theorem 1.4\label{main2} follows from Proposition 4.2\label{main2.2}.~\hfill$\square$

\section{Generalized Quaternion Algebra and Octonion Algebra \\ over a Commutative Ring}\label{quat-oct}\label{section5}
\label{!!!}

Let $K$ be a commutative associative ring with the identity involution and $a$ an invertible element of the ring $R$. We consider the ring $A_1=(K,a)$. Then $A_1$ is a commutative associative ring, since $B=C=I=J=K$, under the notation of Theorem 1.4\label{main2}. It is natural to write elements of the ring $A_1$ in the form $x+yi$, where $x,y$ are elements of the ring $K$, $i=(0,1)$. On the ring $A_1$, an involution is defined by the relation $(x+yi)^*=x-yi$ for any $x,y\in K$. We choose an invertible element $b\in K$. Then $b$ is an invertible symmetrical element of the center of the ring $A_1$ and we can construct the ring $A_2=(A_1,b)$. We consider the $K$-basis of the algebra $A_2$ which is formed by the elements $1=(1,0)$, $i=(i,0)$, $j=(0,1)$ and $k=(0,-i)$. The relations $i^2=a$, $j^2=b$, $ij=-ji=k$, $ik=-ki=aj$, $kj=-jk=bi$ are directly verified. Consequently, the obtained ring is the generalized quaternion algebra $(a,b,K)$ under the notation of \cite{tugan93}. It is well known (and also follows from Theorem 1.3\label{main}) that the ring $A_2$ is associative (e.g., see \cite[Example 7.2.III]{Zhevlakov}). The center of the ring $A_2$ is of the form $K+Ni+Nj+Nk$, where $N=\Ann_K(2)$ (see \cite[Lemma 2(b)]{tugan93}). Let $B,I,J$ be defined  by equations \eqref{centerdefs} for $A=A_2$. It is easy to verify that $B=C=Z(A_2)$, $I=J=N+Ni+Nj+Nk$.

\textbf{Lemma 5.1.}\label{essK}
Under the above notation, the ideal $I$ is an essential ideal in $B$ if and only if $N$ is an essential ideal in $K$.

\textbf{Proof.} Let $I$ be an essential ideal in $B$. If $x\in K\setminus\{0\}$, then there exists an element $y\in B$ such that $xy\in I\setminus\{0\}$. We set
$y=y_1+y_2i+y_3j+y_4k$, where $y_1\in K$ and $y_2,y_3,y_4\in N$. If $xy_1\neq 0$, then $xK\cap N\neq 0$. Otherwise, at least one of the elements $xy_2,xy_3,xy_4$ is not equal to $0$ and each of them belong to the ideal $N$, whence  $xK\cap N\neq 0$ in this case too.

Conversely, if $N$ is an essential ideal in $K$ and $x=x_1+x_2i+x_3j+x_4k\in
I\setminus\{0\}$, then $x_2,x_3,x_4\in N$. If $x_1\neq 0$, then there exists an
element $y\in K$ with $yx_1\in N\setminus\{0\}$. Then $yx\in Bx\cap
I\setminus\{0\}$. If $x_1=0$, then $x=1\cdot x\in Bx\cap I$. Thus,
$I$ is a essential ideal of the ring $B$.~\hfill$\square$

From the above argument, we obtain Proposition 5.2.

\textbf{Proposition 5.2.}\label{csQuat}
The quaternion algebra $((K,a),b)$ is a non-commutative centrally essential ring if and only if $\Ann_K(2)$ is a proper essential ideal of the ring $K$.

Now we consider an arbitrary invertible element $c\in K$ and ring $A_3=(A_2,c)$.
We set $f_1=i$, $f_2=j$, $f_3=k$, $f_4=l=(0,1)$, $f_5=(0,-i)$, $f_6=(0,-j)$, $f_8=(0,-k)$. Some it can be directly verified that the basis $\{1,f_1,f_2,\ldots,f_7\}$ of the $K$-module $A_3$ satisfies the relations from \cite{gen_quat_oct} for basis elements of the generalized octonion algebra ${\mathbb O}(\alpha, \beta, \gamma)$ (for $\alpha=-a, \beta=-b, \gamma=-c$).

Similar to Proposition 5.2\label{csQuat}, we obtain Proposition 5.3.

\textbf{Proposition 5.3.}\label{csOct}
The octonion algebra $(((K,a),b),c)$ is a non-associative centrally essential ring if and only if $\Ann_K(2)$ is a proper essential ideal of the ring $K$.

\textbf{5.4. The completion of the proof of Theorem 1.5.}\\
Let $K=\ZZ_4$. We prove that $R=(((K,1),1),1)$ is a non-associative non-commutative centrally essential ring.\\
Indeed, $\Ann_K(2)=2K$ is an essential proper ideal in $K$. Therefore, the non-commutativity of the ring $((K,1),1)$ (and the non-commutativity of the ring $R$ containing $((K,1),1)$) follows from Proposition 5.2\label{csQuat} and the non-associativity of the ring $R$ follows from Proposition 5.3.\label{csOct}.

We note that the ring $R=(((K,1),1),1)$ is an alternative ring and the ring $(R,1)$ is not even a right-alternative ring, i.e., $(R,1)$ does not satisfy the identity $(x,y,y)=0$ \cite[Exercise 7.2.2]{Zhevlakov}. Thus, there exist alternative non-associative finite centrally essential rings and non-alternative finite centrally essential rings.~\hfill$\square$

\section{Open Questions}\label{section6}

\mbox{$\,$}

\textbf{6.1.} Is it true that there exist left $N$-essential rings which are not right $N$-essential?

\textbf{6.2.} Is it true that there exist commutative $N$-essential (equivalently, centrally essential) non-associative rings?

\textbf{6.3.} Is it true that there exist right-alternative centrally essential or $N$-essential non-alternative rings?

\textbf{6.4.} How can we generalize the obtained results to the case of non-unital rings and the case, where the element $\alpha$ in definition 1.2 not supposed to be invertible?


\end{document}